%% file: main.tex
\documentclass[journal]{IEEEtran}
\usepackage{cite}
\usepackage{amsmath,amssymb,amsfonts}
\usepackage{algorithm}
\usepackage{algorithmic}
\usepackage{graphicx}
\usepackage{textcomp}
\usepackage{subcaption}
\usepackage{hyperref}
\hypersetup{hypertex=true,
            colorlinks=true,
            linkcolor=blue,
            anchorcolor=blue,
            citecolor=blue}
\usepackage{cleveref}
\usepackage[normalem]{ulem}
\usepackage{bm}
\usepackage{multirow}
\usepackage{comment}
    
\usepackage{xcolor}

\begin{document}

\title{Quantum Topology Optimization via Quantum Annealing}
\author{Zisheng Ye, Xiaoping Qian, Wenxiao Pan%
\thanks{Z. Ye, X. Qian and W. Pan are with the Department of Mechanical Engineering, University of Wisconsin-Madison, Madison, WI 53706, USA.}%
\thanks{Corresponding author: Wenxiao Pan (email: wpan9@wisc.edu).}}

\maketitle

\begin{abstract}
We present a quantum annealing-based solution method for topology optimization (TO). In particular, we consider TO in a more general setting, i.e., applied to structures of continuum domains where designs are represented as distributed functions, referred to as continuum TO problems. According to the problem's properties and structure, we formulate appropriate sub-problems that can be solved on an annealing-based quantum computer. The methodology established can effectively tackle continuum TO problems formulated as mixed-integer nonlinear programs. To maintain the resulting sub-problems small enough to be solved on quantum computers currently accessible with small numbers of quits and limited connectivity, we further develop a splitting approach that splits the problem into two parts: the first part can be efficiently solved on classical computers, and the second part with a reduced number of variables is solved on a quantum computer. By such, a practical continuum TO problem of varying scales can be handled on the D-Wave quantum annealer. More specifically, we concern the minimum compliance, a canonical TO problem that seeks an optimal distribution of materials to minimize the compliance with desired material usage. The superior performance of the developed methodology is assessed and compared with the state-of-the-art heuristic classical methods, in terms of both solution quality and computational efficiency. The present work hence provides a promising new avenue of applying quantum computing to practical designs of topology for various applications. 
\end{abstract}


\IEEEpeerreviewmaketitle

\input{Introduction}

\input{Preliminary}

\input{Methodology}

\input{NumericalResults}

\input{Conclusion}

\section*{Acknowledgement}
Z. Y. was supported by the University of Wisconsin - Madison Office of the Vice Chancellor for Research and Graduate Education with funding from the Wisconsin Alumni Research Foundation. W. P. was partially supported by the National Science Foundation under Grant No. CMMI-1761068. X. Q. would like to acknowledge the support of National Science Foundation Grants No. 2219931 and No. 1941206.

\bibliographystyle{IEEEtran}
\bibliography{ref}

\end{document}

%% file: Introduction.tex
\section{Introduction}


Topology optimization (TO) is a computational design method that aims to find optimal distribution of material to improve part performance under governing physical equilibrium equations \cite{bendsoe2003topology}. It originates in structural mechanics \cite{bendsoe1988generating} and has been applied to structures of continuum domains where designs are represented as distributed functions and to discrete truss structures where designs are of finite dimensions and represented via discrete parameters. In this paper, we focus on design problems of distributed parameter systems a.k.a. continuum problems. Continuum TO has received growing attention in a variety of fields, including multi-scale structures \cite{wu2021topology}, fluid mechanics \cite{borrvall2003topology,gersborg2005topology, kreissl2011topology,TOBS_TurbulentFlow_2022}, electromagnetics \cite{kiziltas2003topology, nomura2007structural}, photonics \cite{TO_Photonics_2018,TO_Photonics_2021}, quantum devices \cite{TO_Couplers_2021}, and coupled multi-physics problems \cite{dede2009multiphysics, kreissl2010topology, sun20203d, kumar2020topology}. Topologically optimized designs often exhibit complex free-form shapes, and additive manufacturing can produce parts of complex shapes. The recent wide adoption of additive manufacturing technologies have led to further development of TO for a host of applications \cite{liu2018current}, including imposing additive manufacturability constraints during the optimization process \cite{qian2017undercut, mezzadri2018topology}.

As TO has been expanded from structural mechanics to general multiphysics, it is accompanied by quickly growing complexity and vast searching space, resulting in grand computational challenges \cite{mukherjee2021accelerating,TO_Couplers_2021,liu2018current}. 
Quantum computing (QC) is emerging as a new computing paradigm that could be superior to classical computing for a variety of problems, optimization being among them. 
In literature, only a few attempts have been made to apply quantum algorithms for solving TO problems. Maruo et al. concerned the TO of electromagnetic devices \cite{TO_DA_IEEE2022}. For that, a linear optimization problem was formulated, which can be directly cast into a quadratic unconstrained binary optimization (QUBO) formulation through the Biot-Savart's law. The formulated QUBO problem was then solved through simulated annealing. Their results have indicated that the proposed QUBO formulation can significantly reduce the number of iterations required for the optimization process, compared with the commonly used classical method, normalized Gaussian network \cite{sato2015multimaterial}. Thereby, the overall computational cost can be significantly reduced, because much fewer forward evaluations of the objective functions are invoked, which usually dominate the computational cost. However, the resultant QUBO formulation was not implemented in quantum annealers. Thus, the additional overhead due to quantum embedding and the limitations of quantum hardware currently accessible were not addressed in the analysis of computational cost. Also, the benefit of quantum annealing to converge to the optimal or ground state with larger probability than simulated annealing was not utilized. The work by Sato et al. reported the TO for a simple, discrete truss structure with three edges \cite{QCTO_VAE_2022}. The objective functions of all possible configurations composed by three edges were simultaneously evaluated through quantum entanglement. The optimal configuration with the minimum objective function value was determined by applying the variational quantum algorithms (VQAs) that are suitable for the noisy intermediate-scale quantum (NISQ) devices. The VQAs and forward evaluations of the objective function via finite element analysis (FEA) were integrated together when implemented on a quantum computer, eliminating the need of amplitude estimation for the solution of FEA. In those two works, the proposed quantum algorithms can only be applied to linear TO problems, where the field variables such as magnetic flux density and temperature linearly depend on the design variables. Given that, the TO can be greatly simplified when the objective function is just a linear \cite{QCTO_VAE_2022} or quadratic \cite{sato2015multimaterial} function of the field variables. In addition, the volume or material usage constraint is not strictly imposed by any equality or inequality constraint; as a result, one cannot set a target optimal volume or material usage for the optimization, which is not common in practical TO applications. Further, discrete truss structures only represent a small portion of applications of TO. More generally, TO considers structures of continuum domains where designs are represented as distributed functions.

In the present work, we concern how to apply quantum algorithms to \textit{continuum}, \textit{nonlinear} TO, beyond the discrete truss structures or the linear TO problems tackled in literature \cite{sato2015multimaterial,QCTO_VAE_2022}. In addition, the volume constraint is strictly imposed by an equality or inequality constraint, such that any target optimal volume can be practically achieved. To determine the optimal distribution of materials, continuum TO essentially seeks the numerical solution of a partial differential equation (PDE)-constrained optimization problem. To solve it, the solution domain is first discretized with discrete nodes or elements. The distributed design problem hence involves a large number of binary (0/1) optimization variables, typically one variable per discrete node or element. As the field variables nonlinearly depend on the design variables, it results in a large-scale, mixed integer nonlinear programming (MINLP) problem, subject to equality and/or inequality constraints. This type of problem is intrinsically NP-hard \cite{belotti2013mixed}, namely unable to be solved in polynomial time, and hence is challenging to be handled using classical approaches \cite{andreassen2011efficient,liang2019topology,picelli2021101,munoz2011generalized}.

There are two main implementations of quantum computers: quantum annealers and NISQ processors. Quantum annealers embrace quantum annealing techniques, and the NISQ processors are build with gate-based quantum circuits. Due to the noisy nature of quantum gates and the high overhead cost of noise reduction and error mitigation, the quantum circuit that can be reliably executed on NISQ devices must be limited to short. Thus, hybrid quantum-classical algorithms such as VAQs must be developed for solving optimization on NISQ devices, where only a short parameterized circuit is executed, and classical optimizers are leveraged to train the quantum circuit by optimizing the expectation value over the circuit’s parameters \cite{Liu_LayerVQE_2022}. However, in light of the noise-induced barren plateau \cite{Wang_BarrenPlateaus_NC2021,McClean_BarrenPlateaus_NC2018,Holmes_BarrenPlateaus_PRX2022}, the training process calls for an exponentially scaled computational cost \cite{Wang_BarrenPlateaus_NC2021}, and the classical optimization problems corresponding to the training of parameters are NP-hard \cite{bittel2021training}. Therefore, in the present work, we exploit quantum annealing for solving the TO problem of our interest. Quantum annealing has been demonstrated for its advantages and speedups in solving a variety of combinatorial optimization problems \cite{ajagekar2020quantum, albash2018adiabatic, ajagekar2022hybrid, guillaume2022deep, kelany2022quantum}. For a n-qubit system, it works in the binary discrete space, with the operator defined as a Hamiltonian that can correspond to the objective function of an optimization problem. Quantum annealing can effectively escape local minima by tunneling through the barriers when exploring the objective function's landscape and seeking the global minimum. Thus, it can efficiently solve large-scale nonlinear optimization problems over binary discrete spaces. The final state of the n-qubit system at the end of annealing could be the ground state of the Hamiltonian, with each qubit a classical object standing for the solution to the binary optimization problem. 

The currently accessible quantum computers that are designed to implement quantum annealing are the D-Wave systems \cite{dwave-documentation}. To embed problems on the D-Wave quantum annealer, they need to be formulated in the form of QUBO. However, the optimization problem corresponding to the TO of our interest cannot be directly cast into QUBO formulations. Therefore, we first decompose the original problem into a sequence of mixed integer linear programs (MILPs) via the generalized Benders' decomposition (GBD). Next, in each MILP the continuous variables are represented by a series of auxiliary binary variables, such that each MILP can be mapped onto a binary optimization problem and then through the penalty method transformed into the QUBO formulation. Finally, from the solution of quantum annealing, we obtain the solution of each MILP, and in turn the ultimate optimal material layout after the sequence of MILPs are completely solved. If sufficient, fully connected quantum qubits are available \cite{mcmahon2016fully}, quantum advantages in computational efficiency can be straightforwardly obtained and demonstrated, owing to the efficiency of quantum annealing for solving QUBO. However, currently accessible quantum hardware only offers a small number of qubits with sparse connectivity. To solve practical TO problems in current and near-term quantum annealing devices, we further develop a splitting approach, by which a reduced problem that needs to be solved on a quantum annealer is separated from other sub-problems that can be efficiently handled by classical computers. The resultant, greatly reduced QUBO problem is then implemented on a D-Wave Advantage quantum processing unit (QPU). Both the solution quality and computational cost (including the cost for embedding the problem onto the quantum computer, the time of quantum annealing, and the number of optimization iterations) are analyzed and assessed systematically in a series of problems of variable sizes.

%% file: Preliminary.tex
\section{Preliminary}

\subsection{A Continuum Topology Optimization Problem}
Continuum TO optimizes the material layout $\rho(\mathbf{x})$ within a continuum design domain $\Omega$. A continuum TO problem is usually constrained by a set of partial differential equations (PDEs) subject to boundary conditions (BCs) which describe the physical laws governing the design, as well as the target material usage or volume of material layout. In general, the problem can be stated as:
\begin{equation}
    \begin{aligned}
        \min_{\rho} & \quad \int_{\Omega} F[\mathbf{u}(\rho(\mathbf{x}))] \mathrm{d} \Omega \\
        \text{s.t.} & \quad G_0(\rho(\mathbf{x})) = \dfrac{\int_{\Omega} \rho(\mathbf{x}) \mathrm{d} \Omega}{\int_{\Omega} \mathrm{d} \Omega} - V_T = 0 \\
        & \quad \mathcal{L}(\mathbf{u}(\rho(\mathbf{x}))) + \mathbf{b} = 0 & & \forall \mathbf{x} \in \Omega \\
        & \quad \mathbf{u}(\rho(\mathbf{x})) = \mathbf{u}_\Gamma(\mathbf{x}) & & \forall \mathbf{x} \in \Gamma_D \\
        & \quad \mathbf{n} \cdot \nabla \mathbf{u}(\rho(\mathbf{x})) = \mathbf{h}_\Gamma(\mathbf{x}) & & \forall \mathbf{x} \in \Gamma_N \\
        & \quad G_j(\mathbf{u}(\rho(\mathbf{x})), \rho(\mathbf{x})) = 0, \quad j = 1, \dots m \\
        & \quad H_k(\mathbf{u}(\rho(\mathbf{x})), \rho(\mathbf{x})) \leq 0, \quad k = 1, \dots, n
    \end{aligned}
    \label{eq:continuum_to}
\end{equation}
where $\mathbf{u}$ denotes the field variable defined in the design domain $\Omega$; $\rho$ is the design variable; $\mathcal{L}$ denotes some differential operator; $\mathbf{b}$ is the source term; and $\Gamma = \partial \Omega$ represents the boundary of the design domain with $\Gamma_D$ denoting the boundary where Dirichlet BC is imposed and $\Gamma_N$ the boundary where Neumann BC is enforced. In Eq. \eqref{eq:continuum_to}, $G_0(\rho)=0$ serves to constrain the volume of optimal material layout to the target value $V_T$. The following three lines after $G_0(\rho)=0$ outline the governing PDE and BCs.  $G_j(\mathbf{u}(\rho), \rho)=0$ includes all equality constraints in addition to the volume constraint; and $H_k(\mathbf{u}(\rho), \rho)\leq 0$ comprises all inequality constraints imposed to the optimization problem. In practice, these constraints are usually related to certain manufacture limitations and/or material property requirements. For simplicity, we neglect them in the following discussion. However, they can be easily included in the proposed approach, following the way how we deal with the constraint $G_0(\rho)=0$. The objective function $F$ is convex with respect to the design variable $\rho$, and hence, if the differential operator $\mathcal{L}$ is positive definite, the optimization problem defined in Eq. \eqref{eq:continuum_to} is convex. In case that the TO of interest leads to a non-convex optimization problem, a sequential approximation program can always be employed such that a series of local, convex problems are solved to update the design variable locally in a sequence \cite{nocedal1999numerical}.

To solve the continuum TO problem as in Eq. \eqref{eq:continuum_to}, the design domain $\Omega$ is usually discretized, and the design and field variables are represented in discrete settings. For example, $\Omega$ can be discretized with a uniform mesh, and FEA is used for the numerical discretization of the governing PDEs and BCs. Taking the minimum compliance as an example, the discretized TO problem can be expressed as:
\begin{equation}
    \begin{aligned}
        \min_{\mathbf{u}, \boldsymbol{\rho}} & \quad \mathbf{f}^\intercal \mathbf{u} \\
        \text{s.t.} & \quad \mathbf{K}(\boldsymbol{\rho}) \mathbf{u} = \mathbf{f} \\
        & \quad \sum_{i=1}^{n_{\rho}} \dfrac{\rho_i}{n_{\rho}} = V_T \\
        & \quad \mathbf{u} \in \mathbb{R}^{n_u}, \boldsymbol{\rho} \in \{ 0, 1 \}^{n_{\rho}}
    \end{aligned}
    \label{eq:to_problem}
\end{equation}
where $\mathbf{f}$ is a known external load exerted on the material; the material is subject to linear elasticity; the superscript $^\intercal$ denotes transpose; $\mathbf{u}$ is the discretized displacement field defined on the nodes of the mesh; and $\rho_i$ is the discretized design variable defined on each mesh element $i$. In Eq. \eqref{eq:to_problem},
\begin{equation}
    \mathbf{K}(\boldsymbol{\rho}) = \mathbf{K}_0 + \sum_{i=1}^{n_{\rho}} \rho_i \mathbf{K}_i(E, \nu) \;,
    \label{eq:stiffness_mat}
\end{equation}
where $\mathbf{K}_i \in \mathbb{R}^{{n_u} \times {n_u}}, i = 1, 2, \dots n_{\rho}$ is the predefined element stiffness matrix; $E$ is the Young's modulus; $\nu$ is the Poisson ratio; and $\mathbf{K}_0$ is a symmetric positive definite matrix:
\begin{equation}
    \mathbf{K}_0 = \varepsilon \sum_{i = 1}^{n_{\rho}} \mathbf{K}_i \;,
\end{equation}
with $\varepsilon$ a small number, e.g. $\varepsilon = 10^{-9}$, such that it would not affect the resulting optimal material layout and the values estimated for the field variable on the solid elements with $\rho_i = 1$. The equality constraint $\sum_{i=1}^{n_{\rho}} \dfrac{\rho_i}{n_{\rho}} = V_T$ in Eq. \eqref{eq:to_problem} is referred to as the volume constraint that drives the resulting material layout toward the target volume or material usage. The inclusion of $\mathbf{K}_0$ in $\mathbf{K}$ ensures that for any $\boldsymbol{\rho}$, $\mathbf{K}(\boldsymbol{\rho})$ is a symmetric positive definite matrix. 


The challenge for solving the optimization problem defined in Eq. \eqref{eq:to_problem} lies in the fact that while $\mathbf{u} \in \mathbb{R}^{n_u}$ is a continuous variable, $\boldsymbol{\rho} \in \{ 0, 1 \}^{n_{\rho}}$ is a binary variable, and the filed variable $\mathbf{u}$ nonlinearly depends on the design variable $\boldsymbol{\rho}$, resulting in a MINLP problem. To tackle the challenge, existing classical approaches in TO fall into either of the two categories.  

In the first category, the most widely used method is the Solid Isotropic Material with Penalization (SIMP) method \cite{andreassen2011efficient}. Its strategy is to relax the binary variable into a continuous variable such that the MINLP problem can be converted into a continuous optimization problem. Specifically, the SIMP method approximates the binary design variable with a high-order polynomial function of a continuous variable that smooths out the discontinuity of the binary variable \cite{andreassen2011efficient}. Hence, the density $\rho_i$ continuously varies from 0 to 1 for each element. From the design perspective, such an interpolated representation (with ``gray” density) does not allow easy imposition of design-dependent loading. To recover the desired binary (black/white) representation of material layout in the optimal design, additional numerical treatment is required \cite{sigmund2013topology}, which is not always mathematically justifiable. Also, the relaxation from binary to continuous breaks the convexity of the original problem \eqref{eq:to_problem} and makes the optimization hard to converge \cite{sigmund2013topology}. Further, the resultant continuous optimization problem is not suitable for quantum annealing to solve. Therefore, we do not follow the SIMP method herein, but use its solution as the baseline for comparison with the proposed quantum annealing solution, as discussed in Section \ref{subsec:comparison2other_solver}.

The second category of methods separate the binary and continuous variables into different sub-problems, by iteratively determining one of them with the other fixed in a sub-problem. The representative methods include the Discrete Variable Topology Optimization via Canonical Relaxation Algorithm (DVTOCRA) \cite{liang2019topology}, the Topology Optimization of Binary Structures (TOBS) method \cite{picelli2021101}, and the GBD method \cite{munoz2011generalized}. The DVTOCRA method employs a sequential linear/quadratic approximation to separate the binary and continuous variables into different sub-problems \cite{liang2019topology}. The sub-problem associated with the binary variable is a constrained, quadratic integer programming problem. Due to its NP-hardness, the canonical relaxation algorithm is applied, by which new continuous variables are introduced to replace the binary variable through an approximation function (different from that used in the SIMP method), and the binary optimization programming problem is transformed into a convex, continuous optimization problem in terms of the newly introduced continuous variables. Similarly to the SIMP method, the resultant continuous optimization problems are not suitable for quantum annealing to solve. The TOBS method utilizes sequential linear programming to separate the binary and continuous variables into different sub-problems \cite{picelli2021101}. The resultant sub-problem involving the binary variable is a linear integer programming problem. With the uniform mesh discretization and the exclusion of non-volumetric constraints (i.e., $G_j=0$ and $H_k\le 0$ in Eq. \eqref{eq:continuum_to}), the corresponding integer programming problem is not NP-hard. Hence, the TOBS method does not serve as an ideal candidate for investigating the potential advantages of quantum annealing in TO.

The GBD method follows a different route. It decomposes the original problem into a sequence of MILPs. In contrast to the other methods that have no guarantee of convergence, the GBD method provides a deterministic optimality criterion to warrant convergence within a finite number of iterations \cite{belotti2013mixed}. It is also anticipated to converge faster than the SIMP or TOBS method, because the GBD formulation permits to use all the material layouts generated from previous iterations in each new iteration, while other methods like the SIMP or TOBS can only take into account the material layout yield from the last iteration. The GBD method was first introduced to TO by Mu{\~n}oz et al. \cite{munoz2011generalized} for solving discrete TO problems. In the present paper, we extend it to address a continuum TO problem, as discussed in \Cref{subsec:to_GBD}. In the sub-problem related to the binary design variable after decomposition, additional constraints (associated with the PDEs), other than the volumetric constraint, are introduced and can bring NP-hardness to the optimization problem. Such a sub-problem is well suited to be accelerated by quantum annealing, and hence we focus on the GBD method in the present paper to investigate and demonstrate how QC, particularly via quantum annealing, can be leveraged for solving a continuum, nonlinear TO problem as stated in Eq. \eqref{eq:to_problem}.

\subsection{Generalized Benders' Decomposition}
\label{subsec:to_GBD}

The field variable $\mathbf{u}$ and the design variable $\bm{\rho}$ are separated into different sub-problems via decomposition and are iteratively updated until reaching the convergence. 

At the $k$-th iteration, we have a sequence of feasible integer solutions $\bm{\rho}^j$ that satisfy the volume constraint with \textit{any} desired volume $V$, i.e., $\sum_{i = 1}^{n_{\rho}} \rho_i^j = V$ for $\forall{j}$, and have non-singular stiffness matrices $\mathbf{K}(\bm{\rho}^j)$, where $j$ loops over all previous iterations until the current iteration $k$. We first form the sub-problem with respect to $\mathbf{u}$, i.e., the so-called primal problem, which is given by:
\begin{equation}
    \begin{aligned}
        \min_{\mathbf{u}} & \quad \mathbf{f}^\intercal \mathbf{u} \\
        \text{s.t.} & \quad \mathbf{K}(\boldsymbol{\rho}^{k}) \mathbf{u} = \mathbf{f} \\
        & \quad \mathbf{u} \in \mathbb{R}^{n_u}
    \end{aligned}
    \label{eq:to_primal}
\end{equation}
This primal problem can be easily solved for $\mathbf{u}^k$ by using a linear system solver, as:
\begin{equation}
    \mathbf{u}^k = \mathbf{K}^{-1}(\boldsymbol{\rho}^{k}) \mathbf{f} \;.
    \label{eq:to_fem}
\end{equation}
The primal problem is a restriction to the original problem \eqref{eq:to_problem}, and any $\mathbf{f}^\intercal \mathbf{u}^j$ serves as an upper bound to the problem \eqref{eq:to_problem}.  We denote the lowest upper bound as $U$ such that $U = \min_j (\mathbf{f}^\intercal \mathbf{u}^j), j = 1, \dots, k$. 

Next, we derive the sub-problem with respect to $\boldsymbol{\rho}$, referred to as the master problem. To proceed, we first use a Lagrange multiplier $\bm{\lambda}$ to move the equality constraint to the objective function in problem \eqref{eq:to_problem}, such that the resulting problem is only subject to the constraints involving the binary variables, as:
\begin{equation}
    \begin{aligned}
        \min_{\mathbf{u}, \boldsymbol{\rho}} & \quad \mathbf{f}^\intercal \mathbf{u} + \bm{\lambda}^\intercal \mathbf{K}(\bm{\rho}) \mathbf{u} \\
        \text{s.t.} & \quad \sum_{i=1}^{n_{\rho}} \dfrac{\rho_i}{n_{\rho}} = V \\
        & \quad \mathbf{u} \in \mathbb{R}^{n_u}, \boldsymbol{\rho} \in \{ 0, 1 \}^{n_{\rho}}
    \end{aligned}
    \label{eq:OP_Lagrange}
\end{equation}
The Lagrange multiplier $\bm{\lambda}$ satisfies:
\begin{equation}
    \bm{\lambda}^\intercal \mathbf{K} (\bm{\rho}) = -\mathbf{f}^\intercal \;.
\end{equation}
Then, we introduce an auxiliary variable $\eta$ and replace the objective function in Eq. \eqref{eq:OP_Lagrange} with an inequality constraint, as: 
\begin{equation}
    \begin{aligned}
        \min_{\bm{\rho}, \mathbf{u}, \eta} & \quad \eta \\
        \text{s.t.} & \quad \mathbf{f}^\intercal \mathbf{u} + \bm{\lambda}^\intercal \mathbf{K}(\bm{\rho}) \mathbf{u} \leq \eta \\
        & \quad \sum_{i=1}^{n_{\rho}} \dfrac{\rho_i}{n_{\rho}} = V \\
        & \quad \mathbf{u} \in \mathbb{R}^{n_u}, \boldsymbol{\rho} \in \{ 0, 1 \}^{n_{\rho}}
    \end{aligned}
    \label{eq:to_equiv}
\end{equation}
Applying the first-order Taylor expansion at ($\mathbf{u}^k$,  $\bm{\rho}^k$) and due to the convexity of $\mathbf{f}^\intercal \mathbf{u} + \bm{\lambda}^\intercal \mathbf{K}(\bm{\rho}) \mathbf{u}$, the inequality constraint in Eq. \eqref{eq:to_equiv} can be relaxed as:
\begin{equation}
    \mathbf{f}^\intercal \mathbf{u}^k + \sum_{i = 1}^{n_{\rho}}  {(\bm{\lambda}^k)}^\intercal \mathbf{K}_i \mathbf{u}^k (\rho_i - \rho_i^k) \leq \mathbf{f}^\intercal \mathbf{u} \leq \eta
\end{equation}
By doing so, the problem \eqref{eq:to_equiv} can be relaxed with multiple cuts as:
\begin{equation}
    \begin{aligned}
        \min_{\boldsymbol{\rho}, \eta} & \quad \eta \\
        \text{s.t.} & \quad \mathbf{f}^\intercal \mathbf{u}^{j} + \sum_{i=1}^{n_{\rho}} {(\bm{\lambda}^{j})}^\intercal \mathbf{K}_i \mathbf{u}^{j} (\rho_i - \rho^{j}_i) \leq \eta \\
        & \omit \hfill \text{$j = 1, \dots, k$} \\
        & \quad \sum_{i=1}^{n_{\rho}} \dfrac{\rho_i}{n_{\rho}} = V \\
        & \quad \boldsymbol{\rho} \in \{0, 1\}^{n_{\rho}}
    \end{aligned}
    \label{eq:to_gbd_master}
\end{equation}
Eq. \eqref{eq:to_gbd_master} is the derived master problem for the original TO problem \eqref{eq:to_problem}. The term ${(\bm{\lambda}^j)}^\intercal \mathbf{K}_i \mathbf{u}^j$ denotes the so-called sensitivity in TO. Note that $\bm{\lambda}^j = -\mathbf{u}^j$, due to the symmetry of $\mathbf{K}$, and hence the sensitivity is just $-{(\mathbf{u}^j)}^\intercal \mathbf{K}_i \mathbf{u}^j$. The optimal solution of Eq. \eqref{eq:to_gbd_master} is denoted as $\bm{\rho}^{k+1}$ and added into the sequence of $\bm{\rho}^j$. The optimum of Eq. \eqref{eq:to_gbd_master}, denoted as $\eta^k$, serves as the lower bound of the problem \eqref{eq:to_problem}. The iterations continue until the upper and lower bounds meet such that $(U - \eta^k) / U < \xi$, where $\xi$ is a predefined tolerance for convergence.  

\begin{figure}[htp]
    \centering
    \begin{subfigure}[t]{.45\linewidth}
        \centering
        \includegraphics[width=.5\textwidth]{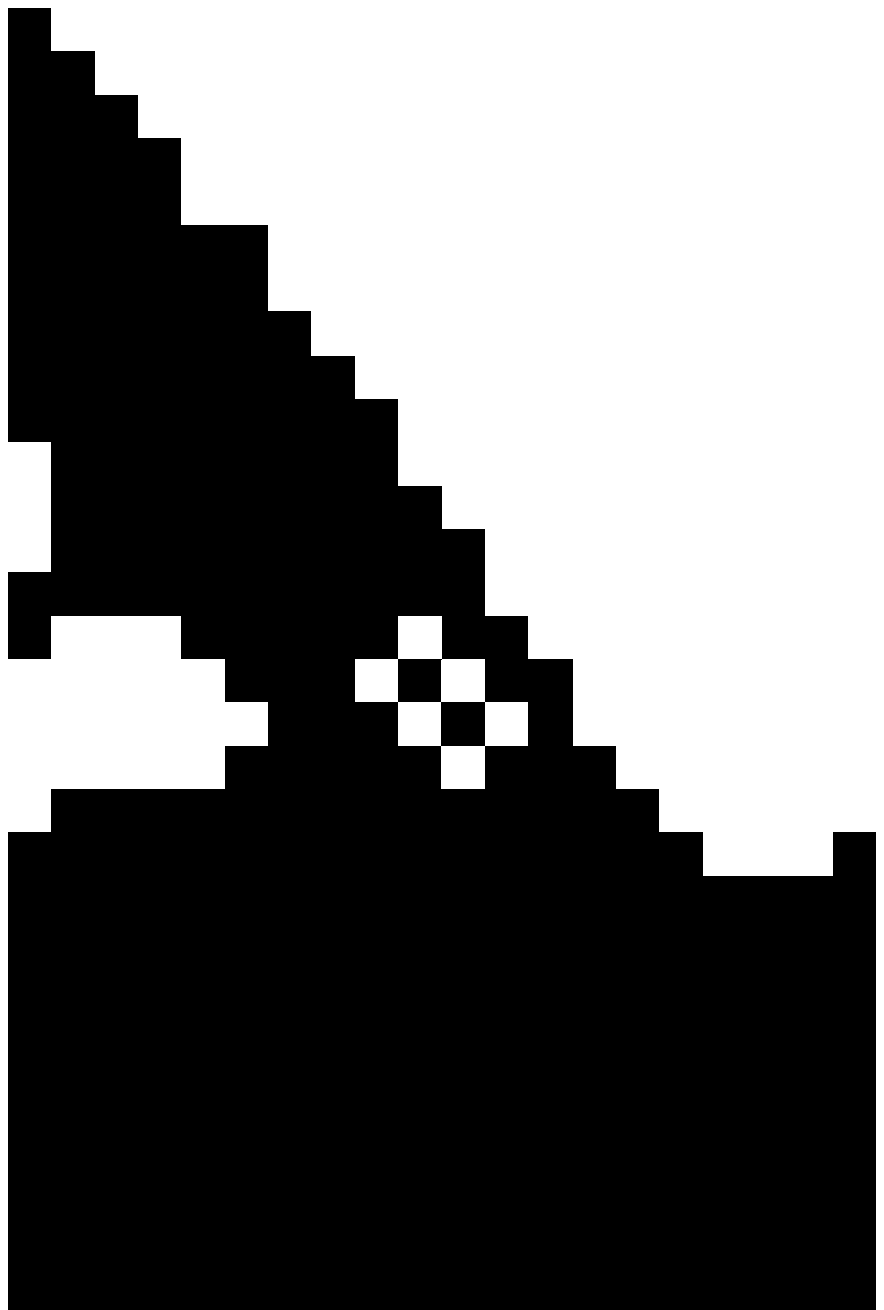}
        \caption{Material layout with ``checkerboard" artifact}
    \end{subfigure}
    \quad \quad
    \begin{subfigure}[t]{.45\linewidth}
        \centering
        \includegraphics[width=.5\textwidth]{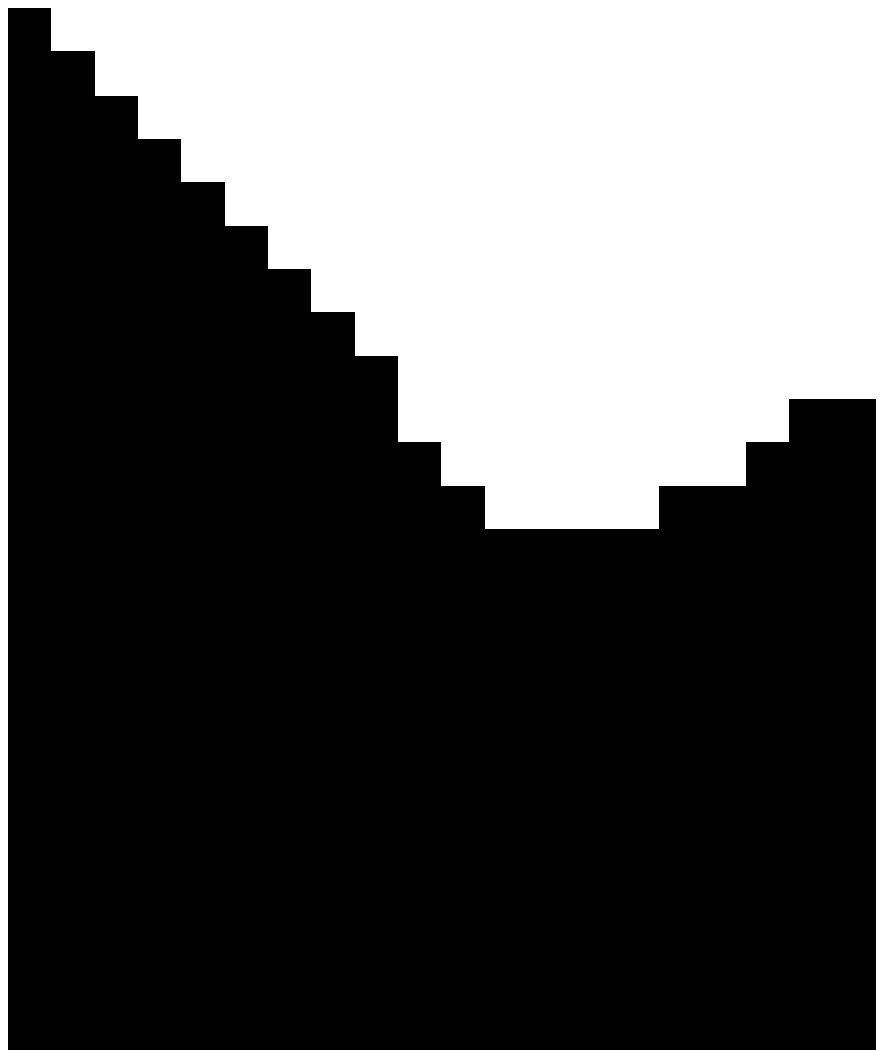}
        \caption{Material layout after filtering}
    \end{subfigure}
    \caption{The effect of filtering on eliminating the checkerboard artifact. The two sub-figures report the resulting optimal material layouts with the same problem setup, using the same optimization procedure, and corresponding to the same region of the solution domain. (a) The optimal material layout yield without filtering exhibits a checkerboard pattern. (b) After applying filtering, the checkerboard artifact is effectively eliminated.}
    \label{fig:checkerboard_illustration}
\end{figure}

Since the continuum TO is essentially a PDE-constrained optimization problem (as in \eqref{eq:continuum_to}), the solution quality also depends on the accuracy of the numerical approximation of the differential operator (e.g., gradient) in the PDE, which becomes especially challenging at the interface between solid ($\rho_i=1$) and void ($\rho_i=0$) elements. The FEA with uniform meshing can lead to the so-called ``checkerboard" artifact, as illustrated in Fig. \ref{fig:checkerboard_illustration}. 
To remedy for that, filtering is required for the sensitivity, following the literature \cite{picelli2021101,sun2022sensitivity}, as: 
\begin{equation}
    \widetilde{w}_i^j =
    \left\{
    \begin{aligned}
        \dfrac{\sum_{l \in \mathcal{N}_i^r} h_{i, l}^r \mathbf{u}^j \mathbf{K}_l \mathbf{u}^j}{\sum_{l \in \mathcal{N}_i^r} h_{i, l}^r}, \quad \rho_i^j = 1 \\
        \varepsilon \dfrac{\sum_{l \in \mathcal{N}_i^r} h_{i, l}^r \mathbf{u}^j \mathbf{K}_l \mathbf{u}^j}{\sum_{l \in \mathcal{N}_i^r} h_{i, l}^r}, \quad \rho_i^j = 0
    \end{aligned}
    \right.
    \label{eq:sensitivity_filtering}
\end{equation}
where $h_{i, l}^r = \max(0, r - \|\mathbf{x}_i - \mathbf{x}_l\|_2)$.
Therefore, the master problem with filtering is rewritten as:
\begin{equation}
    \begin{aligned}
        \min_{\boldsymbol{\rho}, \eta} & \quad \eta \\
        \text{s.t.} & \quad \mathbf{f}^\intercal \mathbf{u}^{j} - \sum_{i=1}^{n_{\rho}} \widetilde{w}_i^j (\rho_i - \rho^{j}_i) \leq \eta \\
        & \omit \hfill \text{$j = 1, \dots, k$} \\
        & \quad \sum_{i=1}^{n_{\rho}} \dfrac{\rho_i}{n_{\rho}} = V \\
        & \quad \boldsymbol{\rho} \in \{0, 1\}^{n_{\rho}}
    \end{aligned}
    \label{eq:to_gbd_master_filter}
\end{equation}

As the iterations proceed, the number of inequality constraints ($j = 1, \dots, k$) included in the master problem \eqref{eq:to_gbd_master_filter} increases. However, those inequality constraints are not independent, and we do not need to include them all when solving the master problem. To accelerate the solution process, we further introduce the Pareto optimal cuts \cite{munoz2011generalized}, as defined by:
\begin{equation}
    \mathcal{P}(k) = \{ j | \mathbf{f}^\intercal \mathbf{u}^j \leq \mathbf{f}^\intercal \mathbf{u}^k, \forall j = 1, \dots k \} \;,
    \label{eq:pareto_cuts}
\end{equation}
which is based on the Pareto dominance relationship and selects only the previous cuts with the objective function values no greater than that of the new cut generated in the current iteration. Such selection of optimal cuts has been numerically proven to improve the rate of convergence to the optimal solution \cite{munoz2011generalized}. Thus, the master problem further with the Pareto optimal cuts can be written as:
\begin{equation}
    \begin{aligned}
        \min_{\boldsymbol{\rho}, \eta} & \quad \eta \\
        \text{s.t.} & \quad \mathbf{f}^\intercal \mathbf{u}^{j} - \sum_{i=1}^{n_{\rho}} \widetilde{w}_i^j (\rho_i - \rho^{j}_i) \leq \eta, \quad \forall j \in \mathcal{P}(k) \\
        & \quad \sum_{i=1}^{n_{\rho}} \dfrac{\rho_i}{n_{\rho}} = V \\
        & \quad \boldsymbol{\rho} \in \{0, 1\}^{n_{\rho}}
    \end{aligned}
    \label{eq:to_master_pareto_cuts}
\end{equation}
If $|P(k)| = 1$, the master problem in Eq. \eqref{eq:to_master_pareto_cuts} can be simplified as:
\begin{equation}
    \begin{aligned}
        \min_{\boldsymbol{\rho}} & \quad \mathbf{f}^\intercal \mathbf{u}^{j} - \sum_{i=1}^{n_{\rho}} \widetilde{w}_i^j (\rho_i - \rho^{j}_i), \quad j \in \mathcal{P}(k) \\
        \text{s.t.} & \quad \sum_{i=1}^{n_{\rho}} \dfrac{\rho_i}{n_{\rho}} = V \\
        & \quad \boldsymbol{\rho} \in \{0, 1\}^{n_{\rho}}
    \end{aligned}
    \label{eq:to_min_master}
\end{equation}
with $\eta^k = \mathbf{f}^\intercal \mathbf{u}^{j} - \sum_{i=1}^{n_{\rho}} \widetilde{w}_i^j (\rho_i^{k+1} - \rho^{j}_i)$.

In summary, for any given volume constraint $V$ and the initial feasible integer solution $\bm{\rho}^1$, the iterative procedure based on the GBD method to solve the TO problem \eqref{eq:to_problem} is outlined in \Cref{algorithm:to_GBD_sub}.
\begin{algorithm}
    \caption{\textsc{ToGbdSub}($V$, $\bm{\rho}^1$)}
    \label{algorithm:to_GBD_sub}
    \textbf{Input}: Volume fraction $V$, initial material layout $\bm{\rho}^1$
    
    \textbf{Output}: Optimal material layout $\boldsymbol{\rho}^*$
    
    \textbf{Utility}: Find an optimal material layout with a given volume fraction $V$ and initial material layout $\bm{\rho}^1$
    
    \begin{algorithmic}[1]
        \FOR{$k=1, \dots$}
            \STATE Initialize the upper bound $U$ as $U = +\infty$
            \STATE Employ the linear system solver to obtain $\mathbf{u}^{k} = \mathbf{K}^{-1}(\boldsymbol{\rho}^{k}) \mathbf{f}$
            \IF {$\mathbf{f}^\intercal \mathbf{u}^{k} < U$}
                \STATE $U \gets \mathbf{f}^\intercal \mathbf{u}^{k}$
                \STATE $\boldsymbol{\rho}^* \gets \boldsymbol{\rho}^{k}$
            \ENDIF
            \STATE Generate the Pareto optimal cuts $\mathcal{P}(k)$ according to \eqref{eq:pareto_cuts}
            \IF {$| \mathcal{P}(k) | = 1$}
                \STATE Solve the master problem \eqref{eq:to_min_master} for the minimizer $\eta^{k}$ and $\boldsymbol{\rho}^{k+1}$
            \ELSE
                \STATE Solve the master problem \eqref{eq:to_master_pareto_cuts} for the minimizer $\eta^{k}$ and $\boldsymbol{\rho}^{k+1}$
            \ENDIF
            \IF {$(U - \eta^k) / U < \xi$}
                \STATE \textbf{break}
            \ENDIF
        \ENDFOR
    \end{algorithmic}
    \hspace*{\algorithmicindent} \textbf{Return}: $\boldsymbol{\rho}^*$
\end{algorithm}


In the following, we explain how to obtain $\bm{\rho}^1$ to initiate the above iterative procedure, which is a feasible binary solution satisfying the volume constraint $\sum_{i=1}^{n_{\rho}} \dfrac{\rho_i}{n_{\rho}} = V$ and has a non-singular stiffness matrix $\mathbf{K}(\bm{\rho}^1)$. First, we consider solving the problem:
\begin{equation}
    \begin{aligned}
        \min_{\boldsymbol{\rho}} & \quad \mathbf{f}^\intercal \mathbf{u}^{0} - \sum_{i=1}^{n_{\rho}} \widetilde{w}_i^0 (\rho_i - \rho^{0}_i) \\
        \text{s.t.} & \quad \sum_{i=1}^{n_{\rho}} \dfrac{\rho_i}{n_{\rho}} = V \\
        & \quad \boldsymbol{\rho} \in \{0, 1\}^{n_{\rho}}
    \end{aligned}
    \label{eq:to_initial_guess}
\end{equation}
where $\mathbf{u}^0 = \mathbf{K}^{-1}(\bm{\rho}^0) \mathbf{f}$;  $\bm{\rho}^0$ can be any material layout that gives a non-singular stiffness matrix $\mathbf{K}(\bm{\rho}^0)$, e.g. with all $\rho_i^0 = 1$, $i = 1, \dots, n_{\rho}$. The problem \eqref{eq:to_initial_guess} is formed by applying the first-order Taylor expansion to the problem \eqref{eq:to_problem} at ($\mathbf{u}^0$, $\bm{\rho}^0$) and with filtering. Hence, it can be regarded as the linear approximation of the original problem. To ensure sufficient accuracy for this approximation, the volume constraint imposed cannot be too far from $\sum_{i=1}^{n_{\rho}} \frac{\rho_i^0}{n_{\rho}}$. However, our ultimate goal is to satisfy the target volume $V_T$, and the value of $V_T$ can generally be far from 1, e.g., $V_T=0.5$. Thus, we adopt an asymptotic process, namely letting the volume constraint gradually approach the target value $V_T$. In particular, we define a sequence of different volume values $V_m = 1 - m \Delta V$ with $m = 1, 2, \dots, M$ and $\Delta V = (1 - V_T)/M$. Limiting $\Delta V$ sufficiently small, we follow an iterative procedure: in the iteration step $m$, solve the problem \eqref{eq:to_initial_guess} with $\bm{\rho}^0$ and letting $V=V_m$, and denote its solution as $\bm{\rho}^1$; invoke \Cref{algorithm:to_GBD_sub} with the input of $\bm{\rho}^1$ to find the optimal solution with respect to the volume constraint $V=V_m$ and let this optimal solution be the new $\bm{\rho}^0$ for the next iteration step $m+1$; and when $m=1$, $\bm{\rho}^0 = \bm{1}$.
The proposed procedure is summarized in \Cref{algorithm:to_GBD}.
\begin{algorithm}
    \caption{\textsc{ToGbd}($V_T$, $\Delta V$)}
    \label{algorithm:to_GBD}
    \textbf{Input}: Target volume fraction $V_T$, incremental change of volume $\Delta V$, number of iterations $M=(1 - V_T)/\Delta V$
    
    \textbf{Output}: Optimal material layout $\boldsymbol{\rho}^*$
    
    \begin{algorithmic}[1]
        \STATE Initialize the design variables as $\boldsymbol{\rho}^{0} = \mathbf{1}$
    	\FOR{$m=1,\cdots M$}
        	\STATE Employ the linear system solver to obtain $\mathbf{u}^{0} = \mathbf{K}^{-1}(\boldsymbol{\rho}^{0}) \mathbf{f}$
            \STATE Solve the problem \eqref{eq:to_initial_guess} to obtain $\boldsymbol{\rho}^{1}$
            \STATE $\bm{\rho}^* \gets \textsc{ToGbdSub}(V_m, \bm{\rho}^1)$
            \STATE Prepare for the next iteration with $\boldsymbol{\rho}^{0} \gets \boldsymbol{\rho}^*$
        \ENDFOR
    \end{algorithmic}
    \hspace*{\algorithmicindent} \textbf{Return}: $\boldsymbol{\rho}^*$
\end{algorithm}

%% file: Methodology.tex
\section{Quantum Algorithm}

\subsection{From Topology Optimization to QUBO}\label{subsec:QUBO}

As discussed above, the original MINLP problem, as in Eq. \eqref{eq:to_problem}, is relaxed into a series of MILP problems, as in Eq. \eqref{eq:to_master_pareto_cuts}. In this section, we establish the conversion of MILP into QUBO, such that the master problem \eqref{eq:to_master_pareto_cuts} can be solved on quantum annealers.

First, a slack variable $\alpha^j$ is introduced into each cut in the problem \eqref{eq:to_master_pareto_cuts} to transform inequality constraints into equality constraints, as:
\begin{equation}
    \begin{aligned}
        \min_{\boldsymbol{\rho}, \eta, \bm{\alpha}} & \quad \eta \\
        \text{s.t.} & \quad \mathbf{f}^\intercal \mathbf{u}^{j} - \sum_{i=1}^{n_{\rho}} \widetilde{w}_i^j (\rho_i - \rho^{j}_i) + \alpha^j = \eta, \quad \forall j \in \mathcal{P}(k) \\
        & \quad \sum_{i=1}^{n_{\rho}} \dfrac{\rho_i}{n_{\rho}} = V \\
        & \quad \boldsymbol{\rho} \in \{0, 1\}^{n_{\rho}} \\
        & \quad \alpha^j \geq 0, \quad \forall j \in \mathcal{P}(k)
    \end{aligned}
\end{equation}
Next, the continuous variables $\eta$ and $\alpha^j$ are replaced with a series of binary variables. To do so, the upper and lower bounds for $\eta$ and $\alpha^j$ need to be specified first.  
Note the facts that $\eta$ is the lower bound of the problem \eqref{eq:to_problem}, and the GBD has generated an upper bound $U$ for \eqref{eq:to_problem}. Thus, $U$ can be regarded as the upper bound for $\eta$. Further, due to the positive definite property of the differential operator $\mathcal{L}$, the compliance defined as $\mathbf{f}^\intercal \mathbf{u}$ must be non-negative. The term, $\mathbf{f}^\intercal \mathbf{u}^{j} - \sum_{i=1}^{n_{\rho}} \widetilde{w}_i^j (\rho_i - \rho^{j}_i)$, is the linear approximation of the original problem \eqref{eq:to_problem} at $(\mathbf{u}^j, \bm{\rho}^j)$, which is close to $\mathbf{f}^\intercal \mathbf{u}^j$ when $\bm{\rho}$ satisfies the volume constraint. Hence, $\mathbf{f}^\intercal \mathbf{u}^{j} - \sum_{i=1}^{n_{\rho}} \widetilde{w}_i^j (\rho_i - \rho^{j}_i) \ge 0$ holds, from which we obtain: $\eta = \mathbf{f}^\intercal \mathbf{u}^{j} - \sum_{i=1}^{n_{\rho}} \widetilde{w}_i^j (\rho_i - \rho^{j}_i) + \alpha^j \ge \alpha^j \ge 0$. Thus, the lower bounds for both $\eta$ and $\alpha^j$ can be set zero. From $\alpha^j \le \eta \le U$, we can set the upper bound for $\alpha^j$ as $U$ as well.  Given their upper and lower bounds specified, the continuous variables $\eta$ and $\alpha^j$ can be approximated by a series of binary variables, as:
\begin{equation*}
    \begin{aligned}
        \eta \approx & \widetilde{\eta}(\mathbf{e}) = U \left( 1 - \dfrac{1}{2^{n_{\eta}}} \right) e_0 + \sum_{i = 1}^{n_{\eta}} \dfrac{U}{2^i} e_i \\
        \alpha^j \approx & \widetilde{\alpha}^j(\mathbf{a}^j) = U \left( 1 - \dfrac{1}{2^{n_{\alpha^j}}} \right) a_0^j + \sum_{i = 1}^{n_{\alpha^j}} \dfrac{U}{2^i} a^j_{i} \;,
    \end{aligned}
\end{equation*}
where
\begin{equation*}
    \begin{aligned}
        \mathbf{e} \in \{0, 1\}^{n_{\eta} + 1}, \quad \mathbf{a}^j \in \{ 0,1 \}^{n_{\alpha^j} + 1} \;.
    \end{aligned} 
\end{equation*}
The accuracy of this approximation is controlled by $n_{\eta}$ and $n_{\alpha^j}$. Their values can be chosen according to the desired accuracy and/or the number of accessible qubits in quantum annealers. Thus, we obtain a binary programming problem as:
\begin{equation}
    \begin{aligned}
        \min_{\boldsymbol{\rho}, \mathbf{e}, \mathbf{a^j}} & \quad \widetilde{\eta}(\mathbf{e}) \\
        \text{s.t.} & \quad \mathbf{f}^\intercal \mathbf{u}^{j} - \sum_{i=1}^{n_{\rho}} \widetilde{w}_i^j (\rho_i - \rho^{j}_i) + \widetilde{\alpha}^j (\mathbf{a}^j) \\
        & \omit \hfill \text{$= \widetilde{\eta}(\mathbf{e}), \quad \forall j \in \mathcal{P}(k)$} \\
        & \quad \sum_{i=1}^{n_{\rho}} \dfrac{\rho_i}{n_{\rho}} = V \\
        & \quad \boldsymbol{\rho} \in \{0, 1\}^{n_{\rho}} \;.
    \end{aligned}
    \label{eq:to_binary}
\end{equation}
By such, the constraints can be readily moved into the objective function through the penalty method to obtain the following QUBO formulation:
\begin{equation}
    \begin{aligned}
        \widetilde{\eta}(\mathbf{e}) & + \sum_{j \in \mathcal{P}(k)} A^j \left[ W^{j} - \left( \sum_{i = 1}^{n_{\rho}} \widetilde{w}_i^{j} \rho_i \right) + \widetilde{\alpha}^j(\mathbf{a}^j) - \widetilde{\eta}(\mathbf{e}) \right]^2 \\
        & + B \left( \sum_{i=1}^{n_{\rho}} \dfrac{\rho_i}{n_{\rho}} - V \right)^2 \;,
    \end{aligned}
    \label{eq:to_master_problem_reformulated}
\end{equation}
where
\begin{equation*}
    W^{j} = \mathbf{f}^\intercal \mathbf{u}^{j} + \sum_{i = 1}^{n_{\rho}} \widetilde{w}_i^j \rho_i^j \;.
\end{equation*} 
As such, we have cast the MILP problem \eqref{eq:to_master_pareto_cuts} into QUBO. 

Finally, we determine the penalty factors. Theoretically, the penalty factors in \eqref{eq:to_master_problem_reformulated} should be gradually increased until the optimal solution converges. However, due to the limited machine precision of currently accessible quantum computers (e.g. up to $10^{-6}$ \cite{dwave-datasheet}), the penalty factors cannot be arbitrarily large in practice. On the other hand, the penalty factors scaled with $U$, the upper bound of $\eta$, can effectively direct the quantum annealer to find a solution near the optimal. From our numerical experiments, setting the penalty factors equal to the upper bound $U$ leads to a satisfactory performance, as a trade-off between the machine precision attainable in the hardware used and the magnitude of the penalty factors required for accuracy,
 i.e., $A^j = U, \quad B = U$. It is worth to mention that any inexact solution for the problem \eqref{eq:to_master_pareto_cuts}, smaller than the exact solution due to the violation of the constraints, will not lead to an abnormal termination of the GBD iterations in \Cref{algorithm:to_GBD_sub}, but just result in more iteration steps before reaching the convergence such that $(U - \eta^k) / U < \xi$.

\subsection{A Splitting Approach for Problem Reduction}
\label{subsec:heuristic_reduction}

Solving the QUBO problem in \eqref{eq:to_master_problem_reformulated} requires $n_{qubit}^l = n_{\rho} + n_{\eta} + \sum_{j \in \mathcal{P}(k)} n_{\alpha^j} + |\mathcal{P}(k)| + 1$ logical qubits  to embed the problem, which can be a large number. Furthermore, the quadratic terms $[ W^{j} - ( \sum_{i = 1}^{n_{\rho}} \widetilde{w}_i^{j} \rho_i ) + \widetilde{\alpha}^j - \widetilde{\eta} ]^2$ and $( \sum_{i=1}^{n_{\rho}} \frac{\rho_i}{n_{\rho}} - V )^2$ require all-to-all connections, since each evolves all binary variables simultaneously. However, due to the limited connectivity of qubits on the currently available quantum computers \cite{lucas2019hard}, all-to-all connections would require much more physical qubits than logical qubits to embed the QUBO problem, making the hardware accessibility even more challenging. Consequently, a practical TO problem formulated as \eqref{eq:to_master_problem_reformulated} cannot be handled by the current or near-term quantum computers. 

To resolve this issue, we further develop a splitting approach in this section. We note that a binary programming problem like \eqref{eq:to_min_master} or \eqref{eq:to_initial_guess} can be efficiently solved by a classical optimizer, as discussed below.
In fact, by taking the continuous relaxation of the design variable as $\bm{\bar{\rho}}$ ($0 \leq \bar{\rho}_i \leq 1, ~ i=1, \dots, n_\rho$), the binary programming problem like \eqref{eq:to_min_master} or \eqref{eq:to_initial_guess} can be relaxed into a (continuous) linear programming. If we consider the volume constraint $\sum_{i=1}^{n_{\rho}} \bar{\rho}_i = n_{\rho} V$ as a cut onto the $n_\rho$-dimensional cube of the relaxed design variable $\bm{\bar{\rho}}$, where $n_{\rho} V$ is an integer, the feasible set for the linear programming can be represented by a polyhedron $P$ as:
\begin{equation*}
    P = \left\{ (\bar{\rho}_1, \dots, \bar{\rho}_{n_{\rho}}) \left| \sum_{i=1}^{n_{\rho}} \bar{\rho}_i = n_{\rho} V, \quad  0 \leq \bar{\rho}_i \leq 1, \forall i \right. \right\}\;,
\end{equation*}
where each element of $\bar{\bm{\rho}}$ must be binary, i.e.,  $\bar{\rho}_i =0$ or 1 $\forall i$, at each vertex. Thus, the linear programming's optimal solution, i.e., a vertex of the polyhedron $P$ \cite{bertsimas1997introduction}, must also be the optimal solution of the original binary programming problem. 
Therefore, the problem \eqref{eq:to_min_master} or \eqref{eq:to_initial_guess} can be exactly solved with the complexity of $\mathcal{O}(n_{\rho})$ and hence be handled efficiently on classical computers. In contrast, although the problem \eqref{eq:to_binary} is also a binary programming problem, the additional $|\mathcal{P}(k)|$ constraints or cuts, other than the volume constraint, make it generally impossible to find the feasible set of the relaxed linear programming as a polyhedron with vertices taking binary values, and hence it still suffers from the NP-hardness.

By leveraging this observation, we propose the following procedure to split the problem \eqref{eq:to_master_pareto_cuts} into two parts. The first part consists of several sub-problems like \eqref{eq:to_min_master}, which will be solved on classical computers. The second part is a problem similar to \eqref{eq:to_master_pareto_cuts} but with greatly reduced variables, which will be solved on a quantum computer. By taking only one of the $|\mathcal{P}(k)|$ inequality constraints as well as the volume constraint in \eqref{eq:to_master_pareto_cuts}, we can form $|\mathcal{P}(k)|$ sub-problems like \eqref{eq:to_min_master} as:
\begin{equation}
    \begin{aligned}
        \min_{\boldsymbol{\rho}} & \quad \mathbf{f}^\intercal \mathbf{u}^{j} - \sum_{i=1}^{n_{\rho}} \widetilde{w}_i^j (\rho_i - \rho^{j}_i), \quad \forall j \in \mathcal{P}(k) \\
        \text{s.t.} & \quad \sum_{i=1}^{n_{\rho}} \dfrac{\rho_i}{n_{\rho}} = V \\
        & \quad \boldsymbol{\rho} \in \{0, 1\}^{n_{\rho}} \;.
    \end{aligned}
    \label{eq:sub_problem_to}
\end{equation}
The solution of \eqref{eq:sub_problem_to} is denoted as $\widetilde{\bm{\rho}}^j$. Note that when $|\mathcal{P}(k)|>1$, the material layouts obtained in different iterations (e.g., $j_1, j_2 \in \mathcal{P}(k)$) are only different in finite numbers of nodes, i.e., most elements of the design variables can be the same. We denote the index set of those elements as $\mathcal{I}$ such that
\begin{equation}
    \widetilde{\rho}_i^{j_1} = \widetilde{\rho}_i^{j_2}, \quad \forall i \in \mathcal{I} \subseteq \{ 1, \dots, n_{\rho} \}, \quad \forall j_1, j_2 \in \mathcal{P}(k) \;.
    \label{eq:sub_index_set}
\end{equation}
And the complement set of $\mathcal{I}$ is denoted as $\mathcal{I^C} = \{1, \dots, n_{\rho} \} \setminus \mathcal{I}$, which contains the index of elements that have different values for the design variable in different iterations, i.e., $\widetilde{\rho}_i^{j_1} \neq \widetilde{\rho}_i^{j_2}, \forall i \in \mathcal{I}^C, \exists j_1, j_2 \in \mathcal{P}(k) $. 

Thus, a reduced problem of \eqref{eq:to_master_pareto_cuts} can be formed as:
\begin{equation}
    \begin{aligned}
        \min_{\boldsymbol{\rho}, \eta} & \quad \eta \\
        \text{s.t.} & \quad R^j - \sum_{i \in \mathcal{I^C}} \widetilde{w}_i^j (\rho_i - \rho^{j}_i) \leq \eta, \quad \forall j \in \mathcal{P}(k) \\
        & \quad \sum_{i \in \mathcal{I^C}} \dfrac{\rho_i}{n_{\rho}} = V - \sum_{i \in \mathcal{I}} \dfrac{\widetilde{\rho}_i}{n_{\rho}} \\
        & \quad \boldsymbol{\rho} \in \{0, 1\}^{|\mathcal{I^C}|} \;,
    \end{aligned}
    \label{eq:reduced_to_master}
\end{equation}
where
\begin{equation}
    R^j = \mathbf{f}^\intercal \mathbf{u}^j - \sum_{i \in \mathcal{I}} \widetilde{w}_i^j (\widetilde{\rho}_i^j - \rho_i^j)\;.
\end{equation}
Compared with the original problem \eqref{eq:to_master_pareto_cuts}, while the number of constraints is identical, the number of design variables has been reduced, from $n_{\rho}$ to $|\mathcal{I^C}|$. Although it is difficult to estimate the exact value of $|\mathcal{I^C}|$ in advance, we anticipate $|\mathcal{I^C}|\ll n_{\rho}$ as discussed above and also demonstrated by the numerical experiments in \Cref{subsec:num_exp_direct}. 

By introducing slack variables $\alpha^j$ and approximating the continuous variables $\eta$ and $\alpha^j$ as series of binary variables, as described in \Cref{subsec:QUBO}, the reduced master problem \eqref{eq:reduced_to_master} can be transformed into a binary programming problem as:
\begin{equation}
    \begin{aligned}
        \min_{\boldsymbol{\rho}, \mathbf{e}, \mathbf{a}_j, \mathbf{b}} & \quad \widetilde{\eta}(\mathbf{e}) \\
        \text{s.t.} & \quad R^j - \sum_{i \in \mathcal{I^C}} \widetilde{w}_i^j (\rho_i - \rho^{j}_i) + \widetilde{\alpha}^j(\mathbf{a}^j) \\
        & \omit \hfill \text{$ = \widetilde{\eta}(\mathbf{e}), \quad j \in \mathcal{P}(k)$} \\
        & \quad \sum_{i \in \mathcal{I^C}} \dfrac{\rho_i}{n_{\rho}} = \widetilde{V} \\
        & \quad \boldsymbol{\rho} \in \{0, 1\}^{|\mathcal{I^C}|}, \mathbf{e} \in \{0, 1\}^{n_{\eta} + 1} \\
        & \quad \mathbf{a}^j \in \{ 0,1 \}^{n_{\alpha^j} + 1} \;,
    \end{aligned}
    \label{eq:reduce_to_master_binary}
\end{equation}
where
\begin{equation*}
    \begin{aligned}
        \widetilde{\eta}(\mathbf{e}) = & U \left( 1 - \dfrac{1}{2^{n_{\eta}}} \right) e_0 + \sum_{i = 1}^{n_{\eta}} \dfrac{U}{2^i} e_i \;, \\
        \widetilde{\alpha}^j(\mathbf{a}^j) = & U \left( 1 - \dfrac{1}{2^{n_{\alpha^j}}} \right) a_{j0} + \sum_{i = 1}^{n_{\alpha^j}} \dfrac{U}{2^i} a^j_{i} \;, \\
        \widetilde{V} = & V - \sum_{i \in \mathcal{A}} \dfrac{\widetilde{\rho}_i}{n_{\rho}} \;.
    \end{aligned}
\end{equation*}
And further through the penalty method, it can be written into the QUBO formulation as:
\begin{equation}
    \begin{aligned}
        \widetilde{\eta}(\mathbf{e}) + & \sum_{j \in \mathcal{P}(k)} A^j \left[ R^j - \left( \sum_{i \in \mathcal{I^C}} \widetilde{w}_i^{j} \rho_i \right) + \widetilde{\alpha}_j(\mathbf{a}_j) - \widetilde{\eta}(\mathbf{e}) \right]^2 \\
        + & B \left[ \sum_{i \in \mathcal{I^C}} \dfrac{\rho_i}{n_{\rho}} - \widetilde{V} \right]^2 \;.
    \end{aligned}
    \label{eq:reduce_to_master_problem_reformulated}
\end{equation}
As $|\mathcal{I^C}|\ll n_{\rho}$, the QUBO problem in \eqref{eq:reduce_to_master_problem_reformulated} 
can now be handled by the current or near-term quantum computers with a small number of qubits and limited connectivity. 

The proposed splitting approach is summarized in \Cref{algorithm:to_GBD_Splitting}.
\begin{algorithm}
    \caption{\textsc{ToGbdSubSplitting}($V$, $\bm{\rho}^1$)}
    \label{algorithm:to_GBD_Splitting}
    \textbf{Input}: Volume fraction $V$, initial material layout $\bm{\rho}^1$
    
    \textbf{Output}: Optimal material layout $\boldsymbol{\rho}^*$
    
    \textbf{Utility}: Find an optimal material layout with a given volume fraction $V$ and initial material layout $\bm{\rho}^1$
    
    \begin{algorithmic}[1]
        \FOR{$k=1, \dots$}
            \STATE Initialize the upper bound $U$ as $U = +\infty$
            \STATE Employ the linear solver to obtain $\mathbf{u}^{k} = \mathbf{K}^{-1}(\boldsymbol{\rho}^{k}) \mathbf{f}$
            \IF {$\mathbf{f}^\intercal \mathbf{u}^{k} < U$}
                \STATE $U \gets \mathbf{f}^\intercal \mathbf{u}^{k}$
                \STATE $\boldsymbol{\rho}^* \gets \boldsymbol{\rho}^{k}$
            \ENDIF
            \STATE Generate the Pareto optimal cuts $\mathcal{P}(k)$ according to \eqref{eq:pareto_cuts}
            \IF {$| \mathcal{P}(k) | = 1$}
                \STATE Solve the master problem \eqref{eq:to_min_master} for the minimizer $\eta^{k}$ and $\boldsymbol{\rho}^{k+1}$ on a classical computer
            \ELSE
                \STATE Solve the sub-problem \eqref{eq:sub_problem_to} for the minimizer $\widetilde{\boldsymbol{\rho}}^{j}$ for each element in $\mathcal{P}(k)$ on a classical computer
                \STATE Generate the index sets $\mathcal{I}$ and $\mathcal{I^C}$ according to \eqref{eq:sub_index_set}
                \STATE Solve the reduced QUBO problem \eqref{eq:reduce_to_master_problem_reformulated} for the minimizer $\eta^k$ and $\bm{\rho}^{k+1}$ on a quantum annealer
            \ENDIF
            \IF {$(U - \eta^k) / U < \xi$}
                \STATE \textbf{break}
            \ENDIF
        \ENDFOR
    \end{algorithmic}
    \hspace*{\algorithmicindent} \textbf{Return}: $\boldsymbol{\rho}^*$
\end{algorithm}

%% file: NumericalResults.tex
\section{Numerical results}

\subsection{Solving a Toy Problem for Validation}\label{subsec:validation}
To validate that the method discussed in \Cref{subsec:QUBO} permits a quantum computer to effectively solve the MILP problem like \eqref{eq:to_master_pareto_cuts}, we consider the following mixed-integer programming problem with linear inequality and equality constraints:
\begin{equation}
    \begin{aligned}
        \min_{u, v, w, t} \quad & v + w + t + (u - 2)^2 \\
        \text{s. t. } \quad & v + 2w + t + u \leq 3 \\
        \quad & v + w + t \geq 1 \\
        \quad & v + w = 1 \\
        \quad & v, w, t \in \{ 0, 1 \}, u \in \mathbb{R} \;.
    \end{aligned}
    \label{equ:toy_prob}
\end{equation}
This problem contains one continuous and three binary variables and one equality and two inequality constrains. Its unique optimal solution is known to be: $u = 2, v = 1, w = t = 0$. 

According to the proposed methodology, we first transform this problem into a binary programming problem, and then rewrite it in the formulation of QUBO. By noting that $0 \le u \leq 3$, we can introduce a binary approximation $\widetilde{u}$ for the continuous variable $u$ as:
\begin{equation*}
    \widetilde{u}(\mathbf{e}) = 3\left( 1 - \dfrac1{2^{n_u}} \right) e_0 + \sum_{i = 1}^{n_u} \dfrac3{2^i} e_i \;.
\end{equation*}
Next, we introduce two slack variables $\alpha^1$ and $\alpha^2$ to convert the two inequality constraints into equality constraints.  By noting that $0 \leq \alpha^1 \leq 3$, $0 \leq \alpha^2 \leq 2$, and $\alpha^2 \in \mathbb{Z}$ since the second inequality constraint only contains integers, we can obtain their binary representations $\widetilde{\alpha}^1$ and $\widetilde{\alpha}^2$, respectively, and thereby rewrite the MILP as a binary programming problem:
\begin{equation}
    \begin{aligned}
        \min_{\mathbf{e}, \mathbf{a}^1, \mathbf{a}^2, v, w, t} \quad & v + w + t + (\widetilde{u} - 2)^2 \\
        \text{s. t. } \quad & v + 2w + t + \widetilde{u} + \widetilde{\alpha}^1 = 3 \\
        \quad & -v - w - t + \widetilde{\alpha}^2 = -1 \\
        \quad & v + w = 1 \\
        \quad & v, w, t \in \{ 0, 1 \} \\
        \quad & \mathbf{e} \in \{0, 1\}^{n_u + 1} \\ 
        \quad & \mathbf{a}^1 \in \{0, 1\}^{n_{\alpha^1} + 1}, \mathbf{a}^2 \in \{0, 1\}^2 \;,
    \end{aligned}
    \label{equ:toy_binaryprog}
\end{equation}
where
\begin{equation*}
    \begin{aligned}
        \widetilde{\alpha}^1(\mathbf{a}^1) = & 3\left( 1 - \dfrac1{2^{n_{\alpha^1}}} \right) a^1_0 + \sum_{i = 1}^{n_{\alpha^1}} \dfrac3{2^i} a^1_i \\
        \widetilde{\alpha}^2(\mathbf{a}^2) = & a_0^2 + a_1^2 \;.
    \end{aligned}
\end{equation*}
Finally, through the penalty method, it can be written into the QUBO as:
\begin{equation}
    \begin{aligned}
        & v + w + t + \left[ \widetilde{u} - 2 \right]^2 + A \left[ v + 2w + t + \widetilde{u} - 3 + \widetilde{\alpha}^1 \right]^2 \\
        & + A\left[ -v - w - t + 1 + \widetilde{\alpha}^2 \right]^2 + A\left[ v + w - 1 \right]^2 \;,
    \end{aligned}
    \label{eq:toy_problem_QUBO}
\end{equation}
where the penalty factor is set as $A = 4$ in this problem.  

\begin{table*}[tp]
    \caption{The results for solving the toy problem \eqref{equ:toy_prob}.}
    \label{tab:comparison_toy}
    \centering
    \begin{tabular}{|c|c|c|c|c|c|c|}
        \hline
        ($n_u$, $n_{\alpha^1}$) & $u$ & $v$ & $w$ & $t$ & $\alpha_1$ & $\alpha_2$ \\
        \hline
        (5, 5) & 1.96875 & 1 & 0 & 0 & 0 & 0 \\
        \hline
        (10, 10) & 2.00098 & 1 & 0 & 0 & 0 & 0 \\
        \hline
    \end{tabular}
\end{table*}
The problem \eqref{eq:toy_problem_QUBO} was directly submitted to the D-Wave Advantage QPU. The annealing time was set with 20$\mu s$, and the sampling was repeated for 1000 times. The final state reached at the end of each annealing was recorded, and the state with the lowest energy among all 1000 samplings was regarded as the ground state of the Hamiltonian defined in \eqref{eq:toy_problem_QUBO} and hence the solution to the binary programming problem \eqref{equ:toy_binaryprog}. The continuous variables $u$ and $\alpha^1$ as well as the integer variable $\alpha^2$ of the original problem \eqref{equ:toy_prob} were then recovered from their binary approximations. The final results for the values of different variables are summarized in \Cref{tab:comparison_toy}. To confirm that the results have converged with respect to the parameter settings in the binary approximations, we compare two different groups of values for $n_u$ and $n_{\alpha^1}$. The results agree with the theoretical solution to the original problem \eqref{equ:toy_prob}, and hence the accuracy of the obtained QC-based solution is validated. An additional note about the continuous variables is made as follows: if we treat the original problem \eqref{equ:toy_prob} as a continuous optimization problem and solve it for the continuous variables, with all the binary variables fixed to the solution obtained from solving \eqref{eq:toy_problem_QUBO}, the accuracy can be even better than recovering them from their binary approximations.


\subsection{Solving the Minimal Compliance Problem}
\label{subsec:num_exp_direct}

\begin{figure}[htp]
    \centering
    \includegraphics[width=.45\textwidth]{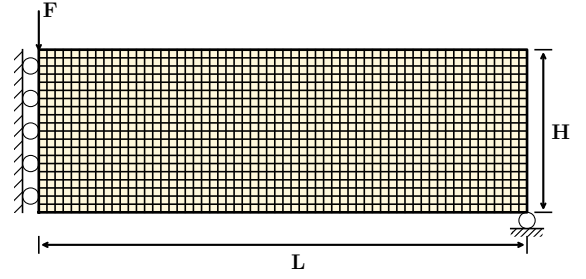}
    \caption{Schematic of the problem setup for the minimum compliance. The yellow region denotes the design domain discretized by a uniform quadrilateral mesh with the resolution of $60 \times 20$, where $L=60$, $H=20$, and the external load is a unit point force, $\mathbf{F}=(0,-1)$, exerted at the top leftmost node.}
    \label{fig:minimum_compliance_illustration}
\end{figure}

In this subsection, we concern a canonical TO problem, the minimum compliance, as formulated in \eqref{eq:to_problem}. The problem setup follows that in \cite{andreassen2011efficient}. More specifically, we consider a rectangular material domain with the external loading of a unit point force $\mathbf{F}$ exerted at the top left corner, as illustrated in \Cref{fig:minimum_compliance_illustration}.  The solution domain, or the design space, is discretized by a uniform quadrilateral mesh with different resolutions. In Eq. \eqref{eq:stiffness_mat}, the Young's modulus is $E = 1.0$, and the Poisson ratio is $\nu = 0.3$. The target volume fraction is $V_T = 0.5$, and the incremental change of volume is $\Delta V = \frac1{24}$.  The convergence criterion tolerance required in \Cref{algorithm:to_GBD_Splitting} is $\xi = 5 \times 10^{-4}$. All values are in reduced units. 

To solve this TO problem, we followed the proposed methodology, consisting of the GBD, conversion of the MILP into QUBO, and the splitting approach. More specifically, \Cref{algorithm:to_GBD} was executed, but with \textsc{ToGbdSub} in Line 5 replaced by \textsc{ToGbdSubSplitting} defined in \Cref{algorithm:to_GBD_Splitting} to integrate the proposed splitting approach. While the binary programming problems \eqref{eq:to_min_master}, \eqref{eq:to_initial_guess} and \eqref{eq:sub_problem_to} were solved by invoking the classical MILP solver Gurobi \cite{gurobi}, the binary programming problem \eqref{eq:reduce_to_master_binary} was solved on the quantum annealer provided by D-Wave. From the validation test in \Cref{subsec:validation} and considering the accuracy of the current quantum devices and the required number of qubits, the parameters for the binary approximations were set as: $n_{\eta} = n_{\alpha^j} = 10$ in \eqref{eq:reduce_to_master_binary}. Once the solution of the design variable $\bm{\rho}^{k+1}$ is determined, the problem \eqref{eq:to_master_pareto_cuts} can be treated as a continuous optimization problem, from which we can solve for the continuous variable $\eta^k$ required by \Cref{algorithm:to_GBD_Splitting}. Alternatively, $\eta^k$ can be recovered from its binary approximation. We chose the former approach herein for its slightly better accuracy.

To implement the problem \eqref{eq:reduce_to_master_binary} on the quantum annealer provided by D-Wave, we examined two different ways. In the first way, we directly embedded the QUBO problem \eqref{eq:reduce_to_master_problem_reformulated} formulated from \eqref{eq:reduce_to_master_binary} in the QPU, which is named as ``\textit{GBD-Splitting-Direct}''. In the second way, we solved the problem \eqref{eq:reduce_to_master_binary} by taking advantage of the hybrid solver provided by D-Wave, particularly the constrained quadratic model (CQM), which is hence referred to as ``\textit{GBD-Splitting-CQM}''. We compare the results in terms of solution quality and wall time, as summarized in \Cref{tab:comparison_direct}, where the discretization resolution is chosen as $60 \times 20$; the parameter for filtering, as described in Eq. \eqref{eq:sensitivity_filtering}, is set as $r = 2$; $T$ denotes the total wall time spent for executing \Cref{algorithm:to_GBD} to find the optimal material layout; $\mathbf{f}^\intercal \mathbf{u}$ denotes the objective function's value obtained corresponding to the optimal material layout; and $N$ denotes the total number of binary programming problems invoked by \Cref{algorithm:to_GBD} and \ref{algorithm:to_GBD_Splitting}. Details about each implementation and our main findings are provided as below.
\begin{table*}[tp]
    \caption{Comparison between different implementations for the reduced binary master problem \eqref{eq:reduce_to_master_binary} on the quantum annealer provided by D-Wave, where $T$ denotes the total wall time spent for executing \Cref{algorithm:to_GBD}.} 
    \label{tab:comparison_direct}
    \centering
    \begin{tabular}{|c|c|c|c|c|c|c|c|c|}
        \hline
        \multirow{2}{*}{\textbf{Resolution}} & \multirow{2}{*}{$n_{\rho}$} &
        \multirow{2}{*}{$r$} & \multicolumn{3}{|c|}{\textbf{GBD-Splitting-Direct}} & \multicolumn{3}{|c|}{\textbf{GBD-Splitting-CQM}}  \\
        \cline{4-9}
        & & & $T$ (s) & $\mathbf{f}^\intercal \mathbf{u}$ & $N$ & $T$ (s) & $\mathbf{f}^\intercal \mathbf{u}$ & $N$ \\
        \hline
        $60 \times 20 $ & 1200 & 2 & 234.12 & 186.3727 & 80 & 111.97 & 178.9243 & 74 \\
        \hline
    \end{tabular}
\end{table*}

In \textit{GBD-Splitting-Direct}, each sampling is set with $20 \mu s$ annealing time, and the sampling is repeated for 1000 times. The final states reached at the end of annealing were recorded, and the state with the lowest energy among 1000 samplings was regarded as the ground state of the Hamiltonian defined in \eqref{eq:reduce_to_master_problem_reformulated} and hence the solution to the binary programming problem \eqref{eq:reduce_to_master_binary}. From our statistics, totally 80 binary programming problems were invoked, 11 among which are the problem \eqref{eq:reduce_to_master_problem_reformulated} and were solved by the quantum annealer. For each, we carefully examined the number of the resultant reduced variables (|$\mathcal{I}^C$|), the number of physical qubits required for embedding the problem \eqref{eq:reduce_to_master_problem_reformulated}, and the solution time spent on the QPU, as summarized in \Cref{tab:direct_embedding_detail}, where the 11 problems are organized according to their sequence of being solved. 
\begin{table*}
    \caption{Statistics about the 11 QUBO problems \eqref{eq:reduce_to_master_problem_reformulated} solved in the implementation of \textit{GBD-Splitting-Direct}, where $T_\text{QUBO}$ denotes the time spent for embedding and annealing on the QPU.}
    \label{tab:direct_embedding_detail}
    \centering
    \begin{tabular}{|c|c|c|c|c|c|c|c|c|c|c|c|}
        \hline
        $|\mathcal{I}^C|$ & 16 & 20 & 34 & 22 & 23 & 27 & 23 & 14 & 32 & 18 & 33 \\
        \hline
        $n_q^r$ & 49 & 53 & 67 & 55 & 56 & 60 & 56 & 47 & 65 & 51 & 66 \\
        \hline
        $n_q^e$ & 198 & 234 & 512 & 260 & 559 & 498 & 432 & 172 & 423 & 212 & 616 \\
        \hline
        $T_\text{QUBO}~(\textrm{s})$ & 6.20 & 7.11 & 16.63 & 8.57 & 29.05 & 21.01 & 32.71 &4.20 & 13.23 & 10.52 & 33.99 \\
        \hline
    \end{tabular}
\end{table*}
With the discretization resolution of $60 \times 20$, the number of design variables is $n_{\rho} = 1200$. From the value of $|\mathcal{I}^C|$ reported in \Cref{tab:direct_embedding_detail}, we can see that the  splitting approach proposed in \Cref{subsec:heuristic_reduction} has greatly reduced the size of the problem to be solved on a quantum computer. Note that the total \textit{logical} qubits needed to represent the problem \eqref{eq:reduce_to_master_problem_reformulated} should also include those required to represent the binary approximations of $\eta$ and $\alpha^j$. Thus, the number of logical qubits required in total is: $n_q^{r} = |\mathcal{I}^C| + n_{\eta} + 1 + |\mathcal{P}(k)| +\sum_{j \in \mathcal{P}(k)} n_{\alpha^j}$. By reducing $|\mathcal{I}^C|$, $n_q^{r}$ can be reduced accordingly. If the \textit{physical} qubits have all-to-all connections, the number of physical qubits required would be consistent with $n_q^{r}$. However, due to the sparse connectivity provided by the current quantum annealing machines, the number of physical qubits required for embedding the problem, denoted as $n_q^e$, is in fact much larger than $n_q^{r}$ \cite{EmbeddingOverhead_PQX2021}, as reported in \Cref{tab:direct_embedding_detail}. Also, with the increase of $n_q^{r}$, $n_q^e$ can grow very fast \cite{EmbeddingOverhead_PQX2021}, although fluctuating due to inhomogeneous connections between qubits. As a result, most of the time spent on the QPU is dominated by the embedding overhead, as indicated by the values of $T_\text{QUBO}$ in \Cref{tab:direct_embedding_detail}, noting that the total annealing time for 1000 repetitions of sampling is only 20 $ms$. The current D-Wave Advantage system permits access to no more than 5000 qubits. To constrain $n_q^e$ not beyond 5000, $|\mathcal{I}^C|$ has to be no more than a few hundreds. Thus, even with the splitting approach, the number of design variables in the original problem must be limited to moderate, which makes solving the minimum compliance problem with finer discretization resolutions challenging.

To tackle this challenge, we pursued the second way of implementation, i.e., \textit{GBD-Splitting-CQM}, where the problem \eqref{eq:reduce_to_master_binary} was submitted, and the CQM, provided by D-Wave, was called for embedding and solving the problem. All default setups were used when employing the CQM
hybrid solver in all the numerical tests. As indicated by the wall time $T$ in \Cref{tab:comparison_direct}, the implementation through CQM greatly saved the entire solving time. The reason for that is the CQM can further reduce the size of the problem embedded on the QPU and in turn save the time spent for embedding, knowing the fact that embedding dominates the time consumed on the quantum devices. In addition, the implementation through CQM results in fewer binary programming problems invoked and a lower value for the objective function. Recalling the discussion about the inexact solution in \Cref{subsec:heuristic_reduction} that any inexact solution to the problem \eqref{eq:to_master_pareto_cuts} can potentially increase the number of GBD iterations, fewer binary programming problem invoked suggest that the solution found by the CQM is closer to the exact solution. Possible reasons for that include: the CQM further reduces the problem embedded on the QPU and makes the resultant optimization problem easier to be accurately solved; in \textit{GBD-Splitting-Direct}, the values chosen for the penalty factors may not be optimal, and the probability of recovering the global optimal solution is highly dependent on the embedding and annealing schedule.

Based on the above comparison and findings, when we scaled up the minimum compliance problem with increasing numbers of design variables (with finer discretization resolutions), we employed the implementation of \textit{GBD-Splitting-CQM}, owing to its better performance in terms of the
solution quality, wall time, and
capacity to embed larger problems. In particular, we considered two different shapes of material domains: $L:H = 3:1$ and $L:H = 2:1$, and the discretization resolution varies from $120\times 40$ to $480\times 240$, leading to the number of design variables ($n_\rho$) varying from 4800 to 115200. The results about the obtained minimum values of the objective function, the associated computing time, and the corresponding optimal material layouts are presented in \Cref{tab:comparison_hybrid} and \Cref{fig:480x240} and \ref{fig:480x160}. To verify the QC-based solutions and to assess the benefits of QC for TO, we further compare with  the state-of-the-art classical solvers, as discussed in the next subsection.

\subsection{Comparison with Classical Solvers}
\label{subsec:comparison2other_solver}
The classical solvers considered include the SIMP method \cite{andreassen2011efficient}, the TOBS method \cite{picelli2021101}, and the DVTOPCRA method \cite{liang2019topology}. For the SIMP method, we followed the implementation in \cite{andreassen2011efficient} with the penalty factor $p = 3$ and the maximum number of iterations set as 1000.  We consider the sensitivity filter controlled by a radius size as in \eqref{eq:sensitivity_filtering} along with Heaviside projection \cite{sigmund2007morphology}.  The convergence criterion is controlled by the change of design variables as $\max_i | \Delta \rho_i | < 0.01$.  For the TOBS method, the setup was the same as provided in \cite{picelli2021101}, except that the MILP solver was changed to Gurobi \cite{gurobi} with default setups to permit a fair comparison with our proposed solution strategy, where Gurobi is used for solving the binary programming problems \eqref{eq:to_min_master}, \eqref{eq:to_initial_guess} and \eqref{eq:sub_problem_to}. For the DVTOCRA method, the implementation exactly followed \cite{liang2020further} for the minimum compliance problem. 

\begin{figure*}[htp]
    \centering
    \begin{subfigure}[t]{.4\textwidth}
        \centering
        \includegraphics[width=.98\textwidth]{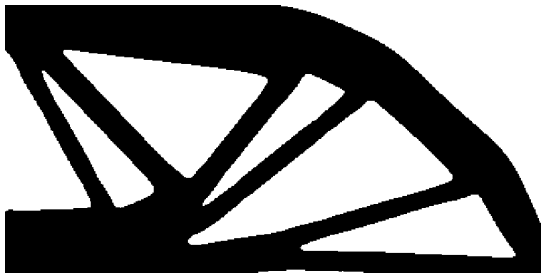}
        \caption{GBD-Splitting-CQM}
    \end{subfigure}
    \quad \quad \quad
    \begin{subfigure}[t]{.4\textwidth}
        \centering
        \includegraphics[width=.98\textwidth]{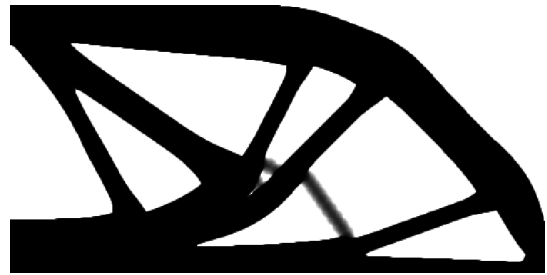}
        \label{subfig:SIMP480x240}
        \caption{SIMP-Heaviside}
    \end{subfigure}
    \par\bigskip
    \begin{subfigure}[t]{.4\textwidth}
        \centering
        \includegraphics[width=.98\textwidth]{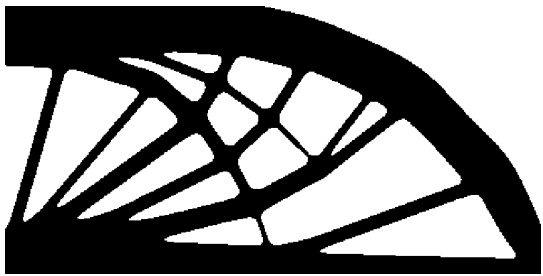}
        \caption{TOBS}
    \end{subfigure}
    \quad \quad \quad
    \begin{subfigure}[t]{.4\textwidth}
        \centering
        \includegraphics[width=.98\textwidth]{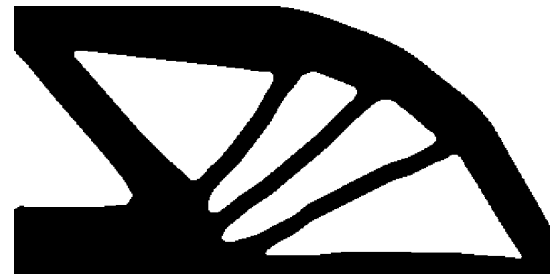}
        \caption{DVTOPCRA}
    \end{subfigure}
    \caption{The resultant optimal material layouts obtained by different methods, with the design space of ($L=60$, $H=30$) and the discretization resolution of $480 \times 240$.}
    \label{fig:480x240}
\end{figure*}

\begin{figure*}[htp]
    \centering
    \begin{subfigure}[t]{.4\textwidth}
        \centering
        \includegraphics[width=.98\textwidth]{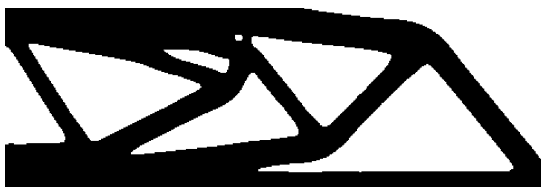}
        \caption{GBD-Splitting-CQM}
    \end{subfigure}
    \quad \quad \quad
    \begin{subfigure}[t]{.4\textwidth}
        \centering
        \includegraphics[width=.98\textwidth]{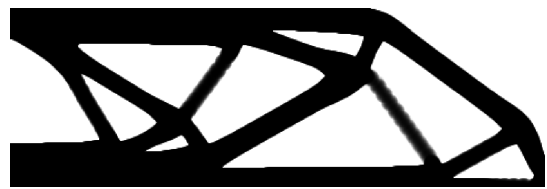}
        \label{subfig:SIMP480x160}
        \caption{SIMP-Heaviside}
    \end{subfigure}
    \par\bigskip
    \begin{subfigure}[t]{.4\textwidth}
        \centering
        \includegraphics[width=.98\textwidth]{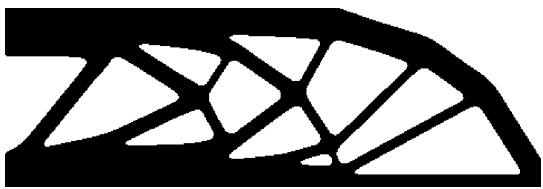}
        \caption{TOBS}
    \end{subfigure}
    \quad \quad \quad
    \begin{subfigure}[t]{.4\textwidth}
        \centering
        \includegraphics[width=.98\textwidth]{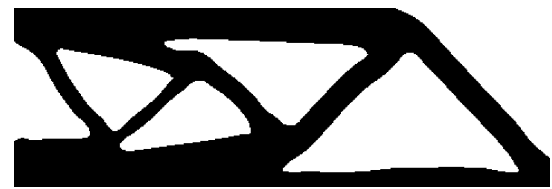}
        \caption{DVTOPCRA}
    \end{subfigure}
    \caption{The resultant optimal material layouts obtained by different methods, with the design space of ($L=60$, $H=20$) and the discretization resolution of $480 \times 160$.}
    \label{fig:480x160}
\end{figure*}

The results obtained by \textit{GBD-Splitting-CQM} and by the three classical methods are compared in \Cref{tab:comparison_hybrid} and \Cref{fig:480x240} and \ref{fig:480x160}, for different shapes of material domains with different discretization resolutions. In \Cref{tab:comparison_hybrid}, $T$ denotes the total wall time spent by each method to find the optimal material layout, including solving all sub-optimization-problems and FEM analysis in \eqref{eq:to_fem}; $\mathbf{f}^\intercal \mathbf{u}$ denotes the objective function's value obtained corresponding to the optimal material layout. For all the three classical methods, $N$ denotes the number of iterations; in each iteration, one FEM analysis \eqref{eq:to_fem} is performed, and one sub-optimization-problem is solved. In \textit{GBD-Splitting-CQM}, $N$ represents the total number of binary programming problems invoked by \Cref{algorithm:to_GBD} and \ref{algorithm:to_GBD_Splitting}; one FEM analysis \eqref{eq:to_fem} is performed prior to invoking each binary programming problem, as stipulated by \Cref{algorithm:to_GBD}. Therefore, comparing the values of $N$ is essentially comparing among all methods how many times the FEM analysis is conducted, which dominates the computational cost in large-scale TO problems.
\begin{table*}[tp]
    \caption{Comparison between different methods for solving the minimum compliance problem.}
    \label{tab:comparison_hybrid}
    \centering
    \begin{tabular}{|c|c|c|c|c|c|c|c|c|c|c|c|c|c|c|c|c|c|}
        \hline
        \multirow{2}{*}{\textbf{Resolution}} & \multirow{2}{*}{$n_{\rho}$} & \multirow{2}{*}{$r$} & \multicolumn{3}{|c|}{\textbf{GBD-Splitting-CQM}} & \multicolumn{3}{|c|}{\textbf{SIMP-Heaviside}} & \multicolumn{3}{|c|}{\textbf{TOBS}} & \multicolumn{3}{|c|}{\textbf{DVTOPCRA}} \\
        \cline{4-15}
        & & &  $T$ (s) & $\mathbf{f}^\intercal \mathbf{u}$ & $N$ & $T$ (s) & $\mathbf{f}^\intercal \mathbf{u}$ & $N$ & $T$ (s) & $\mathbf{f}^\intercal \mathbf{u}$ & $N$ & $T$ (s) & $\mathbf{f}^\intercal \mathbf{u}$ & $N$ \\
        \hline
        $120 \times 40$ & 4800 & 2 & 134.18 & 183.6182 & 74 & 19.84 & 188.1759 & 433 & 6.42 & 187.4646 & 103 & 7.24 & 188.1150 & 139 \\
        \hline
        $120 \times 60$ & 7200 & 2 & 61.51 & 75.1399 & 77 & 37.58 & 77.3887 & 450 & 7.80 & 76.0099 & 83 & 11.03 & 76.7379 & 133 \\
        \hline
        $240 \times 80$ & 19200 & 4 & 123.34 & 183.4272 & 73 & 132.44 & 188.8870 & 530 & 31.82 & 185.7086 & 116 & 39.82 & 191.8400 & 180 \\
        \hline
        $240 \times 120$ & 28800 & 4 & 88.21 & 78.2722 & 89 & 255.84 & 78.8675 & 665 & 39.00 & 77.2155 & 81 & 51.50 & 79.6199 & 145 \\
        \hline
        $480 \times 160$ & 76800 & 8 & 193.90 & 185.1250 & 89 & 1861.68 & 190.7774 & 1000 & 198.14 & 185.1229 & 134 & 259.85 & 191.6287 & 241 \\
        \hline
        $480 \times 240$ & 115200 & 8 & 261.34 & 79.3158 & 90 & 2929.59 & 81.1982 & 1000 & 259.45 & 79.1765 & 95 & 522.30 & 81.8421 & 329 \\
        \hline
    \end{tabular}
\end{table*}

By comparison, the proposed solution strategy exhibits the following advantages. First, it turns out requiring the least number of iterations to reach convergence in almost all the test cases. This can be related to the multi-cuts generated by the GBD iterations and the guarantee of finite iterations by the convexity of the problem. Whereas, the other classical methods solve non-convex sub-problems and in turn have no guarantee for convergence within finite iterations. As the discretization resolution increases, the SIMP method with Heaviside projection cannot meet the convergence criterion within the set maximum number of iterations, and the other two methods, TOBS and DVTOPCRA, tend to require more iterations to reach convergence, as shown in \Cref{tab:comparison_hybrid}. Second, the minimum values of the objective function found by \textit{GBD-Splitting-CQM} are generally lower than those by the other methods, especially SIMP and DVTOPCRA. Lastly, in terms of the wall time, \textit{GBD-Splitting-CQM} exhibits only a slow growth in computing time and hence a good scaling performance, in contrast to the other three methods. Thus, for solving larger-scale problems, e.g., with the discretization resolutions of $480 \times 160$ and $480 \times 240$, \textit{GBD-Splitting-CQM} is more efficient than the SIMP-Heaviside and DVTOPCRA methods and comparable with the TOBS method. Several factors can contribute to that.  First of all, the splitting approach proposed in \Cref{subsec:heuristic_reduction} restricts the growth of the size of the problem to be embedded and solved on a quantum computer and in turn inhibits the growth of time required for embedding.  Next, as the discretization resolution increases, the FEM analysis can be more expensive and in turn dominate the computational cost. As a result, the fewer FEM analysis conducted, the less computing time required. Thus, for higher discretization resolutions,  e.g., $480 \times 160$ and $480 \times 240$, \textit{GBD-Splitting-CQM} is more efficient, as it calls for the least numbers of iterations and hence the least times of FEM analysis. Finally, comparing with the TOBS method, the problem \eqref{eq:to_min_master} contains less constraints than the MILP problem invoked in the TOBS method and in turn costs less solution time in each iteration. 

The resultant optimal material layouts obtained by different methods are shown and compared in \Cref{fig:480x240} and \ref{fig:480x160}, where two different shapes of the design space with the finest discretization resolution are presented. Unlike the nonconvergent SIMP designs with gray regions as shown in Figure \ref{fig:480x240}(b) and \ref{fig:480x160}(b), the \textit{GBD-Splitting-CQM} designs have sharp contrast with clear 0/1 in density. We hence demonstrate the advantage of \textit{GBD-Splitting-CQM} in generating optimal material layouts with clear boundaries, compared against the designs by the method like SIMP that relaxes the binary design variable into a continuous variable. 
In addition, when compared with the designs yield from the TOBS method, the designs optimized by \textit{GBD-Splitting-CQM} also exhibit better minimal length control. Furthermore, the layouts generated by the DVTOPCRA method display crooked internal structures, which can be unfavorable with respect to yield strength.



%% file: Conclusion.tex
\section{Conclusion and Discussion}
We have established a quantum annealing-based solution method for solving continuum, nonlinear TO problems. Its key ingredients include the GBD, conversion of MILP into QUBO, a splitting approach, as well as the implementation on the quantum computer through the hybrid solver, CQM, provided by D-Wave. Through solving a canonical TO problem, the minimum compliance, we have assessed its accuracy and efficiency. By comparing with the state-of-the-art classical methods commonly used in the field of TO, we have demonstrated the advantages of our proposed QC-based methodology, with respect to both solution quality and computational efficiency. Considering the run-time penalty for embedding due to the hardware limitations of current quantum annealing machines \cite{EmbeddingOverhead_PQX2021} and given their future improvements in both the number of qubits and the long-range couplers between qubits, the advantage of the proposed QC-based solution method in terms of computing time can be remarkable. We have also shown its computational efficiency for TO in the sense that the number of FEA runs has been reduced significantly when compared with the classical SIMP method. Future work would exploit this computational efficiency advantages for larger-scale TO problems where FEA cost dominates the optimization run.

This work also presents a new application paradigm for QC and expands the application horizon of QC to include TO, enabling more efficient designs of topology for broad applications. 

With increasing complexity, the continuum TO problems can become more challenging. For example, involving the constraints like $G_j$ and $H_k$ in \eqref{eq:continuum_to} introduces non-volumetric constraints; non-uniform meshes necessitated by discretizing complicated design spaces can lead to non-identical coefficients in the volume constraint associated in the problem \eqref{eq:to_master_pareto_cuts}. The GBD has provided a general enough framework to handle those complexity. Moreover, with non-identical coefficients in the volume constraint, the resultant binary programming problem \eqref{eq:to_min_master} can be NP-hard and cannot be efficiently solved on classical computers. Whereas, the proposed quantum annealing-based solution method, as for the problem \eqref{eq:reduced_to_master}, provides a promising approach for solving NP-hard binary programming problems. In addition, due to the similarities between the problem \eqref{eq:to_min_master} and the sub-problems generated by the TOBS method \cite{picelli2021101}, the proposed solution strategy can also be extended to the TOBS method for dealing with the TO problems subject to complex constraints \cite{sivapuram2018topologyVolumetric}. In light of these perspectives, introducing QC, particularly quantum annealing, to TO may lead to broad and significant impacts. 

The splitting approach developed in the present work has been shown effective for greatly reducing the size of the problem to be embedded in quantum computers, while keeping the solution quality for the optimal material layout. Other approaches, e.g., the multilevel hybrid framework introduced for generic combinatorial optimization \cite{ushijima2021multilevel} and the hierarchical mesh refinement approach \cite{svanberg2007sequential}, can be explored in the future and integrated into the GBD framework to further reduce the cost for solving the problems like \eqref{eq:to_master_pareto_cuts} on quantum computers.

